\documentclass[a4paper,12pt,reqno]{amsart}
\usepackage{amssymb}
\usepackage{amsmath}

\def\i{\,\lrcorner\,}
\def\a{\alpha}

\def\g{\gamma}

\def\la{\langle}
\def\ra{\rangle}
\def\.{\cdot}

\def\n{\nabla}
\def\nb{\bar\nabla}
\def\l{\lambda}

\def\beq{\begin{equation}}
\def\eeq{\end{equation}}
\def\bea{\begin{eqnarray*}}
\def\eea{\end{eqnarray*}}
\def\beaa{\begin{eqnarray}}
\def\eeaa{\end{eqnarray}}
\def\ba{\begin{array}}
\def\ea{\end{array}}

\def\f{\varphi}

\def\o{\omega}

\def\e{\varepsilon}
\def\L{\Lambda}

\def\bp{\begin{proof}}
\def\r{\end{proof}}

\def\Im{{\mathrm {Im}}}
\def\ci{ C^\infty}

\def\CP{{\CM\rm{P}}}

\def \RM{\mathbb{R}}

\def \NM{\mathbb{N}}
\def \ZM{\mathbb{Z}}
\def \CM{\mathbb{C}}

\def\End{{\rm End}}

\def\NK{{\mathcal {NK}}}

\def\d{{\delta}}

\def\Ric{\mathrm{Ric}}

\def\be{\begin{equation}}

\def\ee{\end{equation}}

\def\tr{\mathrm{tr}}
\def\Aut{\mathrm{Aut }}
\def\Cas{\mathrm{Cas}}

\def\RP{\mathbb{R}\mathrm{P}}

\def\Hom{\mathrm{Hom}}
\def\Sym{\mathrm{Sym}}

\def\so{\mathfrak{so}}
\def\su{\mathfrak{su}}

\def\gg{\mathfrak{g}}
\def\tt{\mathfrak{t}}
\def\kk{\mathfrak{k}}
\def\mm{\mathfrak{m}}
\def\uu{\mathfrak{u}}
\def\u{\mathfrak{u}}
\def\pp{\mathfrak{p}}
\def\gl{\mathfrak{gl}}
\def\nn{\mathfrak{n}}
\def\SU{\mathrm{SU}}
\def\SO{\mathrm{SO}}
\def\Sp{\mathrm{Sp}}
\def\U{\mathrm{U}}

\def\ad{\mathrm{ad}}

\def\span{\mathrm{span}}
\def\psp{\Psi ^+}
\def\psm{\Psi ^-}

\def\Sym{\mathrm{Sym}}
\def\scal{\mathrm{scal}}
\def\Ad{\mathrm{Ad}}
\def\Id{\mathrm{id}}

\newtheorem{ede}{Definition}[section]

\newtheorem{epr}[ede]{Proposition}

\newtheorem{ath}[ede]{Theorem}

\newtheorem{elem}[ede]{Lemma}

\newtheorem{ere}[ede]{Remark}

\newtheorem{ecor}[ede]{Corollary}

\title{The Hermitian Laplace Operator on Nearly K\"ahler Manifolds}

\author{Andrei Moroianu and Uwe Semmelmann}

\address{Andrei Moroianu \\ CMLS\\ {\'E}cole Polytechnique \\ UMR 7640 du CNRS
\\ 91128 Palaiseau \\ France}
\email{am@math.polytechnique.fr}

\address{Uwe Semmelmann\\ Mathematisches Institut, Universit{\"a}t zu
K{\"o}ln\\
Weyertal 86-90 D-50931 K{\"o}ln, Germany}
\email{uwe.semmelmann@math.uni-koeln.de}

\thanks{This work was supported by the French-German cooperation
  project Procope no. 17825PG}

\date{\today}

\begin{document}

\begin{abstract}
The moduli space $\NK$ of infinitesimal deformations of a nearly
K\"ahler structure on a compact 6-dimensional manifold is described by
a certain eigenspace of the Laplace operator acting on co-closed
primitive $(1,1)$ forms (c.f. \cite{dnk}). 
Using the Hermitian Laplace operator and some representation theory,
we compute the space $\NK$ on all 6-dimensional homogeneous nearly K\"ahler
manifolds. It turns out that the nearly K\"ahler structure is rigid
except for the flag manifold $F(1,2)=\SU_3/T^2$, which carries an
$8$-dimensional moduli space of infinitesimal nearly K\"ahler
deformations, modeled on the Lie algebra $\su_3$ of the isometry group.
\bigskip


\end{abstract}

\maketitle

\section{Introduction}

Nearly K\"ahler manifolds were introduced in the 70's by A. Gray
\cite{gray} in the context of weak holonomy. More recently,
6-dimensional nearly K\"ahler manifolds turned out to be related to a
multitude of topics among which we mention: Spin manifolds with
Killing spinors (Grunewald), $\SU_3$-structures, geometries with
torsion (Cleyton, Swann), stable forms (Hitchin),
or super-symmetric models in theoretical physics (Friedrich, Ivanov).  

Up to now, the only sources of compact examples are the
naturally reductive 3-symmetric
spaces, classified by Gray and Wolf \cite{gray-wolf}, and the twistor
spaces over positive quaternion-K\"ahler manifolds, 
equipped with the non-integrable almost complex structure. Based on
previous work by R. Cleyton and A. Swann \cite{cs}, 
P.-A. Nagy has shown in 2002 that every simply connected nearly K\"ahler
manifold is a Riemannian product of factors which are either of one of
these two types, or 6-dimensional \cite{nagy}. Moreover, J.-B. Butruille has
shown \cite{jbb} that every homogeneous 6-dimensional nearly K\"ahler
manifold is a 3-symmetric space $G/K$, more precisely isometric with
$S^6=G_2/\SU_3$, $S^3\times 
S^3= \SU_2\times\SU_2\times \SU_2/\SU_2$, $\CP^3=\SO_5/\U_2\times
S^1$ or $F(1,2)=\SU_3/T^2$, all endowed with the metric defined by the
Killing form of $G$.

A method of finding new examples is to take some homogeneous nearly
K\"ahler manifold and try to deform its structure.
In \cite{dnk} we have studied the deformation problem for
6-dimensional nearly K\"ahler manifolds $(M^6,g)$ and proved that if $M$ is
compact, and has normalized 
scalar curvature $\scal_g=30$, then the space $\NK$ of 
infinitesimal deformations of the nearly K\"ahler structure is
isomorphic to the eigenspace for the 
eigenvalue $12$ of the restriction of the Laplace operator $\Delta ^g$ to
the space of {\em co-closed} primitive $(1,1)$-forms $\L^{(1,1)}_0M$. 

It is thus natural to investigate the Laplace
operator on the known 3-symmetric examples (besides the sphere $S^6$,
whose space of nearly K\"ahler structures is well-understood, and
isomorphic to $\SO_7/G_2\cong
\RP^7$, see \cite{fr} or \cite[Prop. 7.2]{jbb}). Recall that the spectrum
of the Laplace operator on symmetric spaces can be computed in
terms of Casimir eigenvalues using the Peter-Weyl formalism. It turns
out that a similar method can be applied in order to compute the
spectrum of a modified Laplace operator $\bar\Delta$ (called the
Hermitian Laplace operator) on 3-symmetric spaces. 
This operator is $\SU_3$-equivariant and 
coincides with the usual Laplace operator on co-closed primitive
$(1,1)$-forms. The space of 
infinitesimal nearly K\"ahler deformations is thus identified with the
space of co-closed forms in
$
\Omega^{(1,1)}_0(12):=\{\a\in \ci(\L^{(1,1)}_0M)\ |\
\bar\Delta\a=12\a\}.
$
Our main result is that the nearly K\"ahler structure is rigid on $S^3\times
S^3$ and $\CP^3$, and that the space of 
infinitesimal nearly K\"ahler deformations of the flag manifold
$F(1,2)$ is eight-dimensional.

The paper is organized as follows. After some preliminaries on nearly
K\"ahler manifolds, we give two general procedures for constructing
elements in	 $\Omega^{(1,1)}_0(12)$ out of Killing vector fields or
eigenfunctions of the Laplace operator for the eigenvalue 12
(Corollary \ref{p1} and Proposition \ref{p2}). We show
that these elements can not be co-closed, thus obtaining an upper
bound for the dimension of the space of infinitesimal nearly K\"ahler
deformations 
(Proposition \ref{direct}). We then 
compute this upper bound explicitly on the  3-symmetric examples and
find that it vanishes for $S^3\times
S^3$ and $\CP^3$, which therefore have no infinitesimal nearly
K\"ahler deformation. This upper bound is equal to $8$ on the flag
manifold $F(1,2)=\SU_3/T^2$ and in the last section we construct an explicit
isomorphism between the Lie algebra of the isometry group $\su_3$ and
the space of infinitesimal nearly K\"ahler deformations on $F(1,2)$.

In addition, our explicit computations (in Section 5) of the spectrum of the
Hermitian Laplace operator on the 3-symmetric spaces, together with
the results in \cite{enk} show that every infinitesimal Einstein
deformation on a 3-symmetric space is automatically an infinitesimal
nearly K\"ahler deformation. 

\medskip

{\it Acknowledgments.} We are grateful to Gregor Weingart for helpful
discussions and in  particular for suggesting the statement of
Lemma~\ref{spectrum}. 

\section{Preliminaries on nearly K\"ahler manifolds}

An almost
Hermitian manifold $(M^{2m},g,J)$ is called {\em nearly K\"ahler} if
\beq\label{1nk}(\n_XJ)(X)=0,\qquad\forall\ X\in TM,\eeq
where $\n$ denotes the Levi-Civita connection of $g$.
The canonical Hermitian connection $\nb$, defined by
\beq\label{der}\nb_XY:=\n_XY-\tfrac12J(\n_XJ)Y,\qquad\forall\
X\in TM,\ \forall\ Y\in C^\infty(M)\eeq
is a $\U_m$ connection on $M$ ({\em i.e.} $\nb g =0$ and $\nb J=0$) with
torsion $\bar T_XY=-J(\n_XJ)Y$. 
A fundamental observation, which goes
back to Gray, is the fact that $\bar\n\bar T=0$ on every nearly
K\"ahler manifold (see \cite{bm}).

We denote the K\"ahler form of $M$ by $\o:=g(J.,.)$. The tensor
$\psp:=\n\o$ is totally skew-symmetric and of type $(3,0)+(0,3)$ by 
(\ref{1nk}). From now on we assume that the dimension of $M$ is $2m=6$ and that
the nearly K\"ahler structure is strict, {\em i.e.} $(M,g,J)$ is not
K\"ahler. It is well-known that $M$ is Einstein in this case. We
will always normalize the scalar curvature of $M$ to $\scal=30$, in
which case we also have $|\psp|^2=4$ point-wise.
The form $\psp$ can be seen as the real
part of a $\nb$-parallel complex volume form $\psp+i\psm$ on $M$,
where $\psm=*\psp$ is the Hodge dual of $\psp$. Thus $M$ carries a 
$\SU_3$ structure whose minimal connection (cf. \cite{cs}) is exactly
$\nb$. Notice that Hitchin has shown that a $\SU_3$ structure
$(\o,\psp,\psm)$ is nearly K\"ahler if and only if the following
exterior system holds:
\beq\label{gray}\begin{cases}d\o=3\psp\\d\psm=-2\o\wedge\o.\end{cases}\eeq

Let $A\in\L^1M\otimes\End M$ denote the tensor
$A_X:=J(\n_XJ)=-\psp_{JX}$, where $\psp_Y$ denotes the endomorphism
associated to $Y\lrcorner\,\psp$ via the metric. Since for every unit vector $X$, $A_X$
defines a complex structure on the 4-dimensional space $X^\perp\cap
(JX)^\perp$, we easily get in a local orthonormal basis $\{e_i\}$ the
formulas
\beq\label{a0}|A_X|^2=2|X|^2,\qquad\forall\ X\in TM.
\eeq
\beq\label{a1}A_{e_i}A_{e_i}(X)=-4X,\qquad\forall\ X\in TM ,
\eeq
where here and henceforth, we use Einstein's summation convention on
repeating subscripts. The following algebraic relations are satisfied
for every $\SU_3$ structure $(\o,\psp)$ on $TM$ (notice that we
identify vectors and 1-forms via the metric):
\beq\label{a10}A_Xe_i\wedge e_i\i\psp=-2X\wedge\o,\qquad\forall\ X\in TM.
\eeq
\beq\label{a3}X\i\psm=-JX\i\psp,\qquad\forall\ X\in TM ,
\eeq
\beq\label{a4}(X\i\psp)\wedge\psp=X\wedge\o ^2,\qquad\forall\ X\in TM.
\eeq
\beq\label{a5}(JX\i\psp)\wedge\o=X\wedge\psp,\qquad\forall\ X\in TM.
\eeq
The Hodge operator satisfies $*^2=(-1)^p$ on $\L^pM$ and moreover
\beq\label{a6}*(X\wedge\psp)=JX\i\psp,\qquad\forall\ X\in TM.
\eeq
\beq\label{a7}*(\phi\wedge\o)=-\phi,\qquad\forall \
\phi\in\L^{(1,1)}_0M.
\eeq
\beq\label{a8}*(JX\wedge\o ^2)=-2X,\qquad\forall\ X\in TM.
\eeq
From now on we assume that $(M,g)$ is compact 6-dimensional
not isometric to  the round sphere $(S^6, can)$. It is well-known that
every Killing vector field $\xi$ on $M$ is an automorphism of the whole
nearly K\"ahler structure (see \cite{dnk}). In particular,
\beq\label{a2}L_\xi\o=0,\qquad L_\xi\psp=0,\qquad L_\xi\psm=0.
\eeq

\noindent
Let now $R$ and $\bar R$ denote the curvature tensors of $\n$ and $\bar\n$.
Then the formula (c.f.~\cite{bfgk})
$$
\begin{array}{rcl}
R_{W X Y Z} &=&\bar R_{W X Y Z} -\tfrac14g(Y,\,W)g(X,\,Z)
+\tfrac14g(X,\,Y)g(Z,\,W) 
\\[1.5ex]
&&+\tfrac34g(Y,\,JW)g(JX,\,Z) -\tfrac34g(Y,\,JX)g(JW,\,Z)
-\tfrac12g(X,\,JW)g(JY,\,Z)
\end{array}
$$
may be rewritten as
$$
R_{XY} = -\,X\wedge Y + R^{CY}_{XY}
$$
and
$$
\bar R_{XY} = -\tfrac{3}{4}\,(X\wedge Y + JX \wedge JY
-\tfrac23\omega(X,Y)J) 
+ R^{CY}_{XY}
$$
where $R^{CY}_{XY}$ is a curvature tensor of Calabi-Yau type.

We will recall the definition of the curvature endomorphism $q(R)$
(c.f.~\cite{dnk}). Let $EM$ 
be the vector bundle associated to the bundle of orthonormal frames via a
representation $\pi: \SO(n) \rightarrow \Aut(E)$. The Levi-Civita connection of
$M$ induces a connection on $EM$, whose curvature satisfies
$R^{EM}_{XY} = \pi_*(R_{XY})= \pi_*(R(X\wedge Y))$, where
we denote with $\pi_*$ the differential of $\pi$ and identify the 
Lie algebra of $SO(n)$, i.e. the skew-symmetric endomorphisms,	with
$\Lambda^2$. In order to keep notations as simple as possible, we
introduce the notation $ \pi_*(A) =A_*$. The curvature endomorphism
$q(R) \in \End(EM)$ is defined as 
\beq\label{qr-def}
q(R) = \tfrac12	(e_i\wedge e_j)_* R(e_i\wedge e_j)_*
\eeq
for any local orthonormal frame $\{e_i\}$. In particular, $q(R)=
\Ric$ on $TM$. 
By the same formula we may define for 
any curvature tensor $S$, or more generally any endomorphism $S$ of
$\Lambda^2TM$,	
a bundle morphism $q(S)$. In any point $q: R \mapsto q(R)$ defines an
equivariant	 
map from the space of algebraic curvature tensors to the space of
endomorphisms of $E$. Since a Calabi-Yau algebraic curvature tensor
has vanishing Ricci curvature,	$q(R^{CY})=0$ holds on $TM$. 
Let $R^0_{XY}$ be defined by $R^0_{XY} = X\wedge Y + JX\wedge JY -
\tfrac23\omega(X,Y)J $.	
Then a direct calculation gives
$$
q(R^0)= \tfrac12\sum (e_i\wedge e_j)_* (e_i\wedge e_j)_* + \tfrac12\sum
(e_i\wedge e_j)_* (Je_i\wedge Je_j)_*  
-\tfrac23\omega_*\omega_*  .
$$
We apply this formula on $TM$. The first summand is exactly the
$\SO(n)$-Casimir, which acts as $-5\Id$. 
The third summand is easily seen to be $\tfrac23 \Id$, whereas the
second summand acts as	
$-\Id$ (c.f. \cite{enk}). Altogether we obtain $q(R^0)=
-\tfrac{16}{3}\Id$, which gives the following expression for $q(\bar
R)$ acting on $TM$:
\beq\label{Rbar}
q(\bar R)\vert_{TM} = 4 \,\Id_{TM}.
\eeq

\section{The Hermitian Laplace operator}

In the next two sections $(M^6,g,J)$ will be a compact nearly K\"ahler
manifold with scalar  
curvature normalized to $\scal_g= 30$.  
We denote as usual by $\Delta$ the Laplace operator $\Delta = d^*d + dd^* =
\nabla^*\nabla +q(R)$ 
on differential forms. We introduce the {\em Hermitian
Laplace 
operator}
\beq\label{Delta-bar-def}
\bar \Delta = \bar\nabla^*\bar\nabla + q(\bar R),
\eeq
which can be defined on any associated bundle $EM$. In \cite{enk}  
we have computed the difference of the operators 
$\Delta$ and $\bar \Delta$ on a primitive $(1,1)$-form $\phi$:
\begin{equation}\label{difference2}
(\Delta	 - \bar \Delta) \phi = (Jd^*\phi)\lrcorner \Psi^+  .
\end{equation}
In particular, $\Delta$ and $\bar\Delta$ coincide on co-closed
primitive $(1,1)$-forms. 
We now compute the difference $\Delta - \bar \Delta$ on 1-forms. Using the 
calculation in \cite{enk} (or directly from
(\ref{Rbar})) we have $q(R) -q(\bar R) = \Id$ on $TM$. It
remains to 
compute the operator $P = \nabla^*\nabla -\bar\nabla^*\bar\nabla$ on
$TM$. A direct calculation using (\ref{a1})
gives for every $1$-form $\theta$
\bea
P(\theta) &=&	 -\tfrac14 A_{e_i} A_{e_i} \theta
-  A_{e_i} \bar\nabla_{e_i}\theta = \theta -	A_{e_i} \bar\nabla_{e_i}\theta
=\theta+\tfrac12A_{e_i} A_{e_i} \theta- A_{e_i}\nabla_{e_i}\theta \\
&=& -\theta - A_{e_i}\nabla_{e_i}\theta.
\eea

In order to compute the last term, we introduce the metric adjoint
$\alpha : \Lambda^2M \rightarrow TM$ of the bundle homomorphism
$X\in 
TM \mapsto X\lrcorner \Psi^+ \in  
\Lambda^2M$. 
It is easy to check that $\alpha(X\lrcorner \Psi^+) = 
2X$ (c.f. \cite{dnk}). Keeping in mind that $A$
is totally skew-symmetric, we compute for an arbitrary vector $X\in
TM$ 
\bea \la A_{e_i}(\nabla_{e_i}\theta),X\ra&=& \la
A_X{e_i},\nabla_{e_i}\theta\ra= 
\la A_X,e_i\wedge\nabla_{e_i}\theta\ra=\la A_X,d\theta\ra\\&=&-\la
\psp_{JX},d\theta\ra= -\la JX,\a(d\theta)\ra=\la J\a(d\theta),X\ra, 
\eea
whence $A_{e_i}(\nabla_{e_i}\theta)=J\a(d\theta)$. 
Summarizing our calculations we have proved the following

\begin{epr}\label{1form}
Let $(M^6,g,J)$ be a nearly K\"ahler manifold with scalar 
curvature normalized to $\scal_g= 30$. Then for any 1-form $\theta$ it
holds that
$$
(\Delta -\bar \Delta) \theta = - J \alpha (d\theta)  .
$$ 
\end{epr}

\noindent
The next result is a formula for the commutator of $J$ and
$\,\alpha\circ d\,$ on $1$-forms.

\begin{elem}\label{magic}
For all $1$-forms $\theta$, the following formula holds:
$$
\alpha(d\theta) = 4J\theta + J \alpha(dJ\theta)  .
$$
\end{elem}
\begin{proof}
Differentiating the identity $\theta\wedge \Psi^+ = J\theta \wedge
\Psi^-$ gives 
$d\theta \wedge \Psi^+ = dJ\theta \wedge \Psi^- + 2 J\theta \wedge
\omega^2$. With respect to the $\SU_3$-invariant decomposition
$\L^2M=\L^{(1,1)}M\oplus\L^{(2,0)+(0,2)}M$, we can write
$d\theta = (d\theta)^{(1,1)} + \tfrac12
\alpha(d\theta)\lrcorner \Psi^+$ and 
$dJ\theta = (dJ\theta)^{(1,1)} + \tfrac12
\alpha(dJ\theta)\lrcorner \Psi^+$.  Since the 
wedge product of forms of type $(1,1)$ and $(3,0)$ vanishes we derive the
equation
$$
\tfrac12 (\alpha (d\theta)\lrcorner \Psi^+)\wedge \Psi^+ 
= \tfrac12 (\alpha(dJ\theta)\lrcorner\Psi^+)\wedge 
\Psi^- + 2J\theta \wedge \omega^2  .
$$
Using (\ref{a4}) and (\ref{a5}) we obtain
$$ 
\tfrac12 \alpha(d\theta)\wedge \omega^2 = \tfrac12 J\alpha(dJ\theta)\wedge \omega^2 +2J\theta\wedge \omega^2  .
$$ 
Taking the Hodge dual of this equation and using (\ref{a8}) gives
$
J\alpha(d\theta) = -\alpha(dJ\theta) - 4 \theta   ,
$
which proves the lemma.
\r

\noindent
Finally we note two interesting consequences of Proposition~\ref{1form}
and Lemma~\ref{magic}.
\begin{ecor}\label{1form-cor}
For any closed 1-form $\theta$ it holds that
$$
 (\Delta - \bar \Delta) \theta	=  0,	 \qquad
 (\Delta - \bar \Delta)J\theta	=  4J \theta.	   
$$
\end{ecor}
\begin{proof}
For a closed 1-form $\theta$ Lemma~\ref{1form} directly implies that $\Delta$
and $\bar \Delta$ coincide on $\theta$. For the second equation we use 
Proposition \ref{1form} together with Lemma~\ref{magic} to conclude
$$
(\Delta -\bar\Delta)J\theta = - J\alpha(dJ\theta) = 4J\theta - \alpha(d\theta) = 4J\theta 
$$
since $\theta$ is closed. This completes the proof of the corollary.
\r


\section{Special $\bar\Delta$-eigenforms on nearly K\"ahler manifolds}

In this section we assume moreover that $(M,g)$ is not isometric to the 
standard sphere $(S^6,can)$.
In the first part of this section we will show how to construct  
$\bar\Delta$-eigenforms on $M$ starting from Killing vector fields.

Let $\xi$ be a non-trivial Killing vector field on $(M,g)$, which in
particular implies
$d^*\xi=0$ and $\Delta \xi = 2\Ric(\xi)=10 \xi$.  As an immediate
consequence of the Cartan formula and (\ref{a2}) we obtain 
\beq\label{jxi}
dJ\xi =L_\xi \o - \xi \lrcorner d\omega = -3 \xi \lrcorner \psi^+
\eeq
so by (\ref{a0}), the square norm of $dJ\xi$ (as a 2-form) is 
\beq\label{norm}|dJ\xi|^2=18|\xi|^2.
\eeq
In \cite{monase} we showed already that the vector field $J\xi$ is 
co-closed if $\xi$ is a Killing vector field and has unit length. However 
it turns out that this also holds more generally.

\begin{epr}\label{jxi2}
Let $\xi$ be a Killing vector field on $M$. Then $d^*J\xi = 0.$
\end{epr}
\begin{proof}
Let $dv$ denote the volume form of $(M,g)$. We start with computing
the $L^2$-norm of $d^*J\xi$.  
$$
\begin{array}{rcl}
\| d^*J\xi \|^2_{L^2} &=& \int_M\la d^*J\xi,d^*J\xi\ra dv =
\int_M[\la\Delta J\xi, J\xi\ra -\la d^*dJ\xi, J\xi\ra] dv \\[1ex]
&=&\int_M [\la \nabla^*\nabla J\xi, J\xi\ra + 5|J\xi|^2 - |dJ\xi|^2]dv\\[1ex]
&=&\int_M [|\nabla J\xi|^2 +5 |\xi|^2 - |dJ\xi|^2]dv=\int_M [|\nabla
J\xi|^2 -13|\xi|^2]dv 
\end{array}
$$
Here we used the well-known Bochner formula for 1-forms, i.e. $\Delta \theta = \nabla^*\nabla \theta + \Ric(\theta)$,
with $\Ric(\theta)=5\theta$ in our case. Next we consider the
decomposition of $\n J\xi$ into its symmetric and skew-symmetric parts
$2\nabla J\xi = dJ\xi + L_{J\xi}g$, 
which together with (\ref{norm}) leads to
\beq\label{nxi}
|\nabla J\xi |^2
=
\tfrac14 (|dJ\xi|^2 + |L_{J\xi}g|^2)
=
9|\xi|^2 + \tfrac14  |L_{J\xi}g|^2.
\eeq
(Recall that the endomorphism square norm of a 2-form is twice its
square norm as a form). In order to compute the last norm, we express
$L_{J\xi}g$ as follows: 
\bea L_{J\xi}g(X,Y)&=&g(\n_XJ\xi,Y)+g(X,\n_Y
J\xi)\\&=&g(J\n_X\xi,Y)+g(X,J\n_Y\xi)
+\psp(X,\xi,Y)+\psp(Y,\xi,X)\\
&=&-g(\n_X\xi,JY)-g(JX,\n_Y \xi)=-d\xi ^{(1,1)}(X,JY), 
\eea
whence 
\beq\label{jx}\|L_{J\xi}g\|^2_{L^2}=2\|d\xi ^{(1,1)}\|^2_{L^2}.
\eeq 
On the other hand, as an application of Lemma \ref{magic} together
with Equation (\ref{jxi}) we get 
$\a(d\xi)=4J\xi+J\a(dJ\xi)=-2J\xi$, so	
\beq\label{dxi}d\xi
^{(2,0)}=-J\xi\i\psp.\eeq
Moreover, $\Delta\xi=10\xi$ since $\xi$ is a
Killing vector field, which yields
$$\|d\xi ^{(1,1)}\|^2_{L^2}=\|d\xi\|^2_{L^2}-\|d\xi
^{(2,0)}\|^2_{L^2}=10\|\xi\|^2_{L^2}-2\|\xi\|^2_{L^2}=8\|\xi\|^2_{L^2}.$$
This last equation, together with (\ref{nxi}) and (\ref{jx}) gives
$\|\n J\xi\|^2_{L^2}=13\|\xi\|^2_{L^2}.$	Substituting this into the
first equation proves that $d^*J\xi$ has vanishing $L^2$-norm and thus
that $J\xi$ is co-closed. 
\r

\begin{epr} \label{del} Let $\xi$ be a Killing vector field on $M$. Then
$$
\Delta \xi=10\xi,\qquad\mbox{and}\qquad \Delta J\xi = 18 J\xi.
$$
In particular, $J\xi$ can never be a Killing vector field.
\end{epr}
\begin{proof}
The first equation holds for every Killing vector field on an Einstein
manifold with $\Ric=5\Id$. 
From (\ref{jxi}) we know $dJ\xi = -3\xi \lrcorner \Psi^+$. Hence
the second assertion follows from:
$$
d^*dJ\xi = -*d*dJ\xi\stackrel{(\ref{a6})}{=}- 3*d(J\xi \wedge
\Psi^+)=9*(\xi \wedge  
\omega^2)\stackrel{(\ref{a8})}{=}18J\xi .
$$
\r
Since the differential $d$ commutes with the Laplace operator 
$\Delta$, every Killing vector field $\xi$ defines two
$\Delta$-eigenforms of degree 2:
$$
\Delta d J\xi = 18 dJ\xi \qquad\mbox{and}\qquad \Delta d\xi = 10 d\xi
$$
As a direct consequence of Proposition \ref{del}, together with
formulas (\ref{jxi}), (\ref{dxi}), and	Proposition
\ref{1form} we get:
\begin{ecor} Every Killing vector field on $M$ satisfies
$$
 \bar \Delta\xi = 12\xi, \qquad \bar \Delta J\xi = 12 J\xi.	
$$
\end{ecor}
Our next goal is to show that the $(1,1)$-part of $d\xi$ is a $\bar\Delta$
-eigenform. 
By (\ref{dxi}) we have
\beq\label{beta}d\xi = \phi - J\xi \lrcorner \Psi^+,
\eeq
for some $(1,1)$-form $\phi$. 
Using Proposition \ref{jxi2}, we can write in a local orthonormal
basis $\{e_i\}$:  
$$\la
d\xi,\o\ra=\tfrac12\la d\xi,e_i\wedge Je_i\ra=\la\n_{e_i}\xi,Je_i\ra=
d^*J\xi=0,$$
thus showing that $\phi$ is primitive. The differential of $\phi$ can
be computed from the Cartan formula:
\beq \label{dphi}d\phi\stackrel{(\ref{beta})}{=} d(J\xi \lrcorner
\Psi^++d\xi)
\stackrel{(\ref{a3})}{=} -d(\xi \lrcorner \Psi^-)=
-L_\xi\Psi^-+\xi\lrcorner
d\Psi^-\stackrel{(\ref{a2})}{=}-2\xi\lrcorner \o ^2= 
- 4 J\xi \wedge \omega.
\eeq
From here we obtain
$$
* d\phi = -4*(J\xi \wedge\o) = 4 \xi \wedge \omega  ,
$$
whence
$$
\begin{array}{rcl}
d*d\phi &=& 4 d\xi \wedge \omega - 12 \xi \wedge \Psi^+ 
\stackrel{(\ref{beta})}{=} 4\phi\wedge \omega - 4(J\xi \lrcorner \Psi^+)\wedge 
\omega -12\xi \wedge \Psi^+\\[1ex]
&\stackrel{(\ref{a5})}{=}& 4\phi \wedge \omega - 16 \xi \wedge \Psi^+.
\end{array}
$$
Using (\ref{a6}) and (\ref{a7}), we thus get
$$
d^*d\phi = -*d*d\phi = 4\phi +16J\xi \lrcorner \Psi^+.
$$
On the other hand,
$$
d^*\phi =
-*d*\phi\stackrel{(\ref{a7})}{=}*d(\phi\wedge\o)
\stackrel{(\ref{dphi})}{=}X(-4J\xi\wedge\o	
^2+3\phi\wedge\psp)
\stackrel{(\ref{a8})}{=}8\xi
$$
and finally
$$
dd^*\phi = 8d\xi = 8\phi - 8J\xi \lrcorner \Psi^+.
$$
The calculations above thus prove the following proposition

\begin{epr}\label{phi}
Let $(M^6,g,J)$ be a compact nearly K\"ahler manifold with scalar 
curvature $\scal_g = 30$, not isometric to the standard sphere. Let
$\xi$ be a Killing vector field on $M$ and let $\phi$ be the
$(1,1)$-part of $d\xi$. Then  
$\phi$ is primitive, i.e. $\phi = (d\xi)^{(1,1)}_0$. Moreover 
$d^*\phi = 8\xi$ and
$
\Delta \phi = 12 \phi + 8J\xi \lrcorner \Psi^+ .
$
\end{epr}

\begin{ecor}\label{p1} The primitive $(1,1)$-form $\f$ satisfies
$$
\bar\Delta \phi = 12 \phi.
$$
\end{ecor}
\begin{proof} From (\ref{difference2}) and the proposition above we
  get
$$\bar\Delta \phi = \Delta\phi- (\Delta-\bar\Delta)\phi=
12 \phi + 8J\xi \lrcorner \Psi^+ - (Jd^*\phi) \lrcorner \Psi^+=12\phi.$$
\r

In the second part of this section we will present another way of
obtaining  primitive $\bar\Delta$-eigenforms of type $(1,1)$, starting from 
eigenfunctions of the Laplace operator. Let $f$ be such an
eigenfunction, i.e.	 
$\Delta f = \lambda f$. We consider the primitive $(1,1)$-form $\eta :=
(dJdf)^{(1,1)}_0$. 

\begin{elem} \label{46}
The form $\eta$ is explicitly given by
$$
\eta = dJdf + 2df\lrcorner \Psi^+ + \tfrac{\lambda}{3}f\omega .
$$
\end{elem}
\begin{proof} According to the decomposition of $\L^2M$ into
  irreducible $\SU_3$-summands, we can write 
$$ dJdf=\eta+\gamma \lrcorner \Psi^++h\o$$
for some vector field $\gamma$ and function $h$.
From Lemma~\ref{magic} we get $2\gamma=\alpha(dJdf)= -4df$. In order
to compute $h$, we write
$$6h\,dv=h\o\wedge\o ^2=dJdf\wedge\o ^2=d(Jdf\wedge\o
^2)\stackrel{(\ref{a8})}{=} 2d*df=2\lambda f\,dv.$$
\r

We will now compute the Laplacian of 
the three summands of $\eta$ separately. First, we have $
\Delta df = \lambda df$ and Corollary~\ref{1form-cor} yields
$\bar\Delta df = \lambda df$. 
Since $\bar\Delta$ commutes with $J$, we also have $\bar\Delta Jdf =
\lambda Jdf$ and from the second equation in Corollary~\ref{1form-cor}
we obtain 
$$
\Delta Jdf = \bar \Delta Jdf + (\Delta-\bar\Delta)Jdf = (\lambda +4)Jdf .
$$
Hence, $dJdf$ is a $\Delta$-eigenform for the eigenvalue $\lambda +
4$. 

\begin{elem}  The co-differential of the $(1,1)$-form $\eta$ is given by
$$
d^*\eta = \left(\tfrac{2\lambda}{3}-4\right) Jdf.
$$
\end{elem}
\begin{proof}
Notice that $d^*(f\omega)=-df \lrcorner \omega$ and  that $d^*Jdf = -*d*Jdf = - \tfrac12 * d (df \wedge \omega^2) = 0 $, since $d\omega^2 = 0$.
Using this we obtain
$$
\begin{array}{rcl}
d^*\eta &=& \Delta Jdf + 2 d^*(df\lrcorner \Psi^+)-\tfrac{\lambda}{3}df\lrcorner \omega 
= (\lambda + 4)Jdf -2*d(df\wedge \Psi^-)-\tfrac{\lambda}{3}Jdf \\[1ex]
&=& (\lambda + 4 - \tfrac{\lambda}{3})Jdf - 4*(df\wedge \omega^2) 
\stackrel{(\ref{a8})}{=}
(\tfrac{2\lambda}{3}-4) Jdf.
\end{array}
$$

\r

\noindent
In order to compute $\Delta$ of the second summand of $\eta$ we need
three additional formulas 

\begin{elem}
$$
\bar\Delta (X \lrcorner \Psi^+) = (\bar\Delta X)\lrcorner \Psi^+.
$$
\end{elem}
\begin{proof}
Recall that $\bar\Delta = \bar\nabla^*\bar\nabla + q(\bar R)$.
Since $\Psi^+$ is $\bar\nabla$-parallel we immediately obtain
$$
\bar\nabla^*\bar\nabla (X\lrcorner \Psi^+) = - \bar\nabla_{e_i}\bar\nabla_{e_i}(X\lrcorner \Psi^+)
= -(\bar\nabla_{e_i}\bar\nabla_{e_i}X)\lrcorner \Psi^+  .
$$
The map $A\mapsto A_*\Psi^+$ is a $\SU_3$-equivariant map from $\Lambda^{2}$ to $\Lambda^{3}$. But
since $\Lambda^3$ does not contain the representation $\Lambda^{(1,1)}_0$ as an irreducible
summand, it follows that $A_*\Psi^+=0$ for any skew-symmetric endomorphism $A$ corresponding 
to some primitive $(1,1)$-form. Hence we conclude
$$
q(\bar R) (X \lrcorner \Psi^+) = \omega_{i *} \bar R (\omega_{i})_*(X
\lrcorner \Psi^+) 
=
(\omega_{i *} \bar R (\omega_{i})_*X) \lrcorner \Psi^+
=(q(\bar R) X) \lrcorner \Psi^+,
$$
where, since the holonomy of $\nb$ is included in $\SU_3$, the sum
goes over some orthonormal basis $\{\omega_{i}\}$ of 
$\L^{(1,1)}_0M$. 
Combining these two formulas we obtain $\bar\Delta (X\lrcorner \Psi^+)
= (\bar \Delta X)\lrcorner \Psi^+$. 
\r

\begin{elem}
$$
(\Delta - \bar\Delta) (df \lrcorner \Psi^+) = 6 (df\lrcorner \Psi^+)- \tfrac{4\lambda}{3}f\omega - 2\eta.
$$
\end{elem}
\begin{proof}
From Proposition 3.4 in \cite{enk} we have
\bea
(\Delta - \bar\Delta) (df \lrcorner \Psi^+)&=&(\n^*\n-\nb^*\nb) (df
\lrcorner \Psi^+)+(q(R)-q(\bar R)) (df \lrcorner \Psi^+)\\
&=&(\n^*\n-\nb^*\nb) (df
\lrcorner \Psi^+)+4df \lrcorner \Psi^+.
\eea
The first part of the right hand side reads
\beq\label{o}(\n^*\n-\nb^*\nb) (df
\lrcorner \Psi^+)=-\tfrac14A_{e_i*}A_{e_i*}df \lrcorner
\Psi^+-A_{e_i*}\nb_{e_i}(df \lrcorner \Psi^+).\eeq
From (\ref{a1}) we get
\bea A_{e_i*}A_{e_i*}df \lrcorner
\Psi^+&=& A_{e_i*}(A_{e_i}e_k\wedge\Psi^+(df,e_k,\.))\\
&=&A_{e_i}A_{e_i}e_k\wedge\Psi^+(df,e_k,\.)+A_{e_i}e_k\wedge
A_{e_i}\Psi^+(df,e_k,\.)\\
&=&-4e_k\wedge e_k\i\psp_{df}+A_{e_i}e_k\wedge
A_{e_i}e_j\Psi^+(df,e_k,e_j)=-8\psp_{df},
\eea
where we used the vanishing of the expression $E=A_{e_i}e_k\wedge
A_{e_i}e_j\Psi^+(df,e_k,e_j)$:
\bea E&=& A_{Je_i}e_k\wedge
A_{Je_i}e_j\Psi^+(df,e_k,e_j)=A_{e_i}Je_k\wedge
A_{e_i}Je_j\Psi^+(df,e_k,e_j)\\
&=&A_{e_i}e_k\wedge
A_{e_i}e_j\Psi^+(df,Je_k,Je_j)=-E.
\eea
It remains to compute the second term in (\ref{o}). We notice that by
Schur's Lemma,
every $\SU_3$-equivariant map from the space of symmetric tensors
$\Sym^2M$ to $TM$ vanishes, so in particular (since $\n df$ is
symmetric), one has $A_{e_i}\n_{e_i}df=0$. We then compute
\bea A_{e_i*}\nb_{e_i}\Psi^+_{df}&=&A_{e_i*}((\nb_{e_i}df)
\lrcorner \Psi^+)=(A_{e_i}\nb_{e_i}df) \lrcorner \Psi^++(\nb_{e_i}df)
\lrcorner A_{e_i*}\psp\\
&\stackrel{(\ref{a10})}{=}&(A_{e_i}\n_{e_i}df) \lrcorner \Psi^+-\tfrac12(A_{e_i}A_{e_i}df) \lrcorner \Psi^+-2(\nb_{e_i}df)
\lrcorner(e_i\wedge\o)\\
&=&2\psp_{df}+2d^*df\o+\la A_{e_i}df,e_i\ra\o+2e_i\wedge J\nb_{e_i}df\\
&=&2\psp_{df}+2\l f\o+2e_i\wedge\nb_{e_i}Jdf
=2\psp_{df}+2\l f\o+2dJdf-e_i\wedge A_{e_i}Jdf\\
&=&2\psp_{df}+2\l f\o+2dJdf+2A_{Jdf}=4\psp_{df}+2\l f\o+2dJdf.
\eea
Plugging back what we obtained into (\ref{o}) yields
$$(\n^*\n-\nb^*\nb) (df
\lrcorner \Psi^+)=-(2\psp_{df}+2\l f\o+2dJdf),$$
which together with Lemma \ref{46} and the first equation prove the
desired formula.
\r
\begin{elem}
$$
\Delta f\omega = (\lambda + 12)f\omega - 2(df\lrcorner \Psi^+) .
$$
\end{elem}
\begin{proof}
Since $d^*(f\omega) = -df \lrcorner \omega = -Jdf $ we have $dd^*(f\omega) = -dJdf$. For the second summand
of $\Delta (f\o)$ we first compute $d (f\omega) = df \wedge \omega +
3f \Psi^+$. Since $d^*\Psi^+=\tfrac13d^*d\omega = 4 \omega$, we get $d^*f\Psi^+= -df\lrcorner
\Psi^+ + fd^*\Psi^+ = -df\lrcorner
\Psi^+ + 4f\omega$. Moreover
$$
\begin{array}{rcl}
d^*(df \wedge \omega) &=& -*d(Jdf \wedge \omega) = - *(dJdf \wedge \omega - 3Jdf\wedge \Psi^+)\\[1ex]
&=&
-*([\eta-2df\lrcorner\Psi^+-\tfrac{\lambda}{3}f\omega]\wedge \omega) + 3 *(Jdf \wedge \Psi^+)\\[1ex]
&=&
\eta +2*((df\lrcorner\Psi^+)\wedge \omega) +\tfrac{2\lambda}{3}f\omega -3df\lrcorner\Psi^+\\[1ex]
&=&
\eta + 2 df\i\Psi^+ +
\tfrac{2\lambda}{3}f\omega  - 3df\lrcorner\Psi^+.
\end{array}
$$
Recalling that $\eta = dJdf + 2df\lrcorner \Psi^+ +
\tfrac{\lambda}{3}f\omega$, we obtain
$$
\Delta f\omega = -dJdf -3df\lrcorner\Psi^+ +12f\omega + \eta -df\lrcorner\Psi^+ +
\tfrac{2\lambda}{3}f\omega = (\lambda + 12)f\omega - 2df\lrcorner\Psi^+.
$$
\r

\noindent
Applying these three lemmas we conclude
$$
\Delta (df \lrcorner \Psi^+) = \bar\Delta(df \lrcorner \Psi^+) +
(\Delta -\bar\Delta)(df \lrcorner \Psi^+) 
=(\lambda + 6) (df \lrcorner \Psi^+) - \tfrac{4\lambda}{3}f\omega - 2\eta
$$
and thus
$$
\begin{array}{rcl}
\Delta \eta &=& (\lambda + 4)dJdf + (2\lambda +12)(df \lrcorner
\Psi^+)-\tfrac{8\lambda}{3}f\omega -4\eta 
+\tfrac{\lambda}{3}(\lambda + 12)f\omega -\tfrac{2\lambda}{3}(df
\lrcorner \Psi^+)\\[1ex] 
&=& \lambda \eta + \left( 4- \tfrac{2\lambda}{3} \right)(df\lrcorner \Psi^+).
\end{array}
$$
Finally we have once again to apply the formula for the difference of
$\Delta$ and $\bar\Delta$ on 
primitive $(1,1)$-forms. We obtain
$$
\bar\Delta \eta = \Delta \eta - Jd^*\eta \lrcorner \Psi^+ 
= \Delta \eta + \left( \tfrac{2\lambda}{3}-4\right)(df\lrcorner \Psi^+)
= \lambda \eta  .
$$
Summarizing our calculations we obtain the following result.

\begin{epr}\label{p2}
Let $f$ be an $\Delta$-eigenfunction with $\Delta f = \lambda f$ Then
the	 primitive $(1,1)$-form	 $\eta := (dJdf)^{(1,1)}_0$ satisfies
$$ 
\bar \Delta\eta = \lambda \eta \qquad \mbox{and} \qquad d^*\eta =
\left(\tfrac{2\lambda}{3} -4\right)Jdf .
$$
\end{epr}

Let $\Omega^0(12) \subset C^\infty(M)$ be the
$\bar\Delta$-eigenspace for the eigenvalue 12  
(notice that $\bar\Delta = \Delta$ on functions) and
let $\Omega^{(1,1)}_0(12)$ denote the space of primitive
$(1,1)$-eigenforms of $\bar\Delta$ corresponding to the eigenvalue 12.
Summarizing Corollary \ref{p1} and Proposition \ref{p2}, we have constructed a
linear mapping 
$$\Phi:i(M)\to\Omega^{(1,1)}_0(12),\qquad \Phi(\xi):=d\xi^{(1,1)}_0$$ 
from the space of  Killing vector fields into $\Omega^{(1,1)}_0(12)$ and a
linear mapping
$$\Psi:\Omega^0(12) \to
\Omega^{(1,1)}_0(12),\qquad\Psi(f):=(dJdf)^{(1,1)}_0.$$  
Let moreover $\NK\subset\Omega^{(1,1)}_0(12)$ denote the space of
nearly K\"ahler 
deformations, which by \cite{dnk} is just the space of co-closed forms
in $\Omega^{(1,1)}_0(12)$.
\begin{epr}\label{direct}
The linear mappings $\Phi$ and $\Psi$  defined above  are injective and 
the sum $\Im(\Phi)+\Im(\Psi)+\NK\subset\Omega^{(1,1)}_0(12)$ is a direct
sum. In particular, 
\beq\label{ineq}\dim(\NK)\le \dim(\Omega^{(1,1)}_0(12))-\dim(i(M))-
\dim(\Omega^0(12)).
\eeq
\end{epr}

\begin{proof}
It is enough to show that if $\xi\in i(M)$, $f\in
\Omega^0(12)$ and $\alpha\in\NK$ satisfy 
\beq\label{e1}d\xi^{(1,1)}_0+(dJdf)^{(1,1)}_0+\alpha=0,\eeq 
then $\xi=0$ and $f=0$. 
We apply $d^*$ to (\ref{e1}). Using Propositions \ref{phi} and
\ref{p2} to express the 
co-differentials of the first two terms we get
\beq\label{8}
8\xi+8Jdf=0.
\eeq 
Since $J\xi$ is co-closed (Proposition \ref{jxi2}), formula (\ref{8})
implies $0=d^*J\xi=d^*df=12f$, i.e. $f=0$. Plugging back into
(\ref{8}) yields $\xi=0$ too.
\r

\section{The homogeneous Laplace operator on reductive homogeneous spaces}

\subsection{The Peter-Weyl formalism}
Let $M=G/K$ be a homogeneous space with compact Lie groups $K\subset G$ and 
let $\pi:K\to\Aut(E)$ be a representation of $K$. We denote by
$EM:=G\times_\pi E$ be the associated vector bundle over $M$.  
The Peter-Weyl theorem and the Frobenius reciprocity yield
the following isomorphism of $G$-representations: 
\beq\label{pw} 
L^2(EM)\cong \bigoplus_{\g\in\hat G} V_\g\otimes \Hom_K(V_\g,E), 
\eeq 
where $\hat G$ is the set of (non-isomorphic) irreducible 
$G$-representations. 
If not otherwise stated we will consider only 
complex representations. Recall that the space of smooth sections
$C^\infty(EM)$ 
can be identified with the space $C^\infty(G;E)^K$ of 
$K$-invariant $E$-valued functions,
i.e. functions $f:G\rightarrow E$ with $f(gk)=\pi(k)^{-1}f(g)$. This space is
made into a $G$-representation by the left-regular representation $\ell$,
defined by by $ (\ell(g)f)(a) = f(g^{-1}a)$. Let $v\in V_\g$ and
$A\in \Hom_K(V_\g,E)$ then the invariant $E$-valued function
corresponding to $v\otimes A$ is defined by $g\mapsto A(g^{-1}v)$. In
particular, each summand in the Hilbert space direct sum (\ref{pw}) is
a subset of $C^\infty(EM)\subset L^2(EM)$.

Let $\gg$ be the Lie algebra of $G$. We denote by $B$ the {\em
Killing form} of $\gg$, $B(X,Y):= \tr(\ad_X\circ\ad_Y)$. The Killing 
form is non-degenerated and negative definite if $G$ is compact and 
semi-simple, which will be the case in all examples below.

If $\pi:G\to\Aut(E)$ is a $G$-representation, the {\em Casimir
  operator} of $(G,\pi)$ acts on $E$ by the formula 
\beq\label{f}
\Cas^G_\pi=\sum (\pi_*X_i) ^2,
\eeq
where $\{X_i\}$ is a $(-B)$-orthonormal basis of $\gg$ and
$\pi_*:\gg\rightarrow \End(E)$ denotes the differential of the
representation $\pi$.

\begin{ere} \label{rem1} Notice that the
Casimir operator is divided by $k$ if one use the scalar product $-kB$
instead of $-B$.\end{ere}

\noindent
If $G$ is simple, the adjoint representation $\ad$ on the
complexification $\gg^\CM$ is irreducible, so, by Schur's Lemma, its
Casimir operator  
acts as a scalar. Taking the trace in (\ref{f}) for $\pi=\ad$ yields 
the useful formula $\Cas^G_{\ad}=-1$. 

Let $V_\g$ be an irreducible $G$-representation of highest weight $\g$.
By Freudenthal's formula the Casimir operator acts on $V_\g$ by 
scalar multiplication with $\|\rho\|^2-\|\rho+\g\|^2$, where $\rho$ 
denotes the half-sum of the positive roots and $\|\.\|$ is the norm 
induced by $-B$ on the dual of the Lie algebra of the maximal torus of $G$.
Notice that these scalars are always non-positive. Indeed
$\;\|\rho\|^2-\|\rho+\g\|^2=-\la\g,\g+2\rho\ra_B$ and  
$\la\g,\rho\ra\ge0$, since $\g$ is a dominant weight, i.e. it is in 
the the closure of the fixed Weyl chamber,  whereas $\rho$ is the 
half-sum of positive  weights and thus by definition has a
non-negative scalar product with $\g$.

\subsection{The homogeneous Laplace operator}
We denote by $\bar\nabla$ the canonical homogeneous connection on 
$M=G/K$. It coincides with the Levi-Civita connection only in the 
case that $G/K$ is a symmetric space. A crucial observation is that
the canonical homogeneous connection coincides with the canonical
Hermitian connection on naturally reductive 3-symmetric spaces (see below).
We define the curvature 
endomorphism $q(\bar R) \in \End(EM)$ as in (\ref{qr-def}) and 
introduce as in (\ref{Delta-bar-def}) the second order operator 
$\bar\Delta_\pi = \bar\nabla^*\bar\nabla + q(\bar R)$ acting on 
sections of the associated bundle $EM:=G\times_\pi E$.

\begin{elem} \label{qr}
Let $G$ be a compact semi-simple Lie group, $K\subset G$ a compact subgroup,
and let $M=G/K$ the naturally reductive homogeneous space equipped with the 
Riemannian metric induced by $-B$. For every $K$-representation $\pi$
on $E$, let $EM:=G\times_\pi E$
be the associated vector bundle over $M$. Then the endomorphism 
$q(\bar R) $ acts fibre-wise on $EM$ as $\;q(\bar R) = -\Cas^K_\pi$. Moreover
the differential operator $\bar\Delta$ acts on the space of sections of 
$EM$, considered as $G$-representation via the left-regular
representation,  as $\bar \Delta  = -\Cas^G_\ell$.
\end{elem}
\begin{proof} 
Consider the $\Ad(K)$-invariant 
decomposition
$\gg = \kk \oplus \pp$. For any vector $X\in\gg$ we write
$X=X^\kk+X^\pp$, with $X^\kk \in \kk$ and $X^\pp \in \pp$. The canonical 
homogeneous connection is the left-invariant connection in the principal 
$K$-fibre bundle $G\rightarrow G/K$ corresponding to the projection 
$X\mapsto X^\kk $. It follows that one can do for the canonical homogeneous
connection on $G/K$ the same identifications as for the Levi Civita
connection on Riemannian symmetric spaces. 

In particular, the covariant  derivative of a section $\phi \in \Gamma(EM)$
with respect to some $X\in \pp$ translates into the derivative
$X(\hat\phi)$ of the the corresponding function $\hat\phi\in
C^\infty(G;E)^K$, which is minus the differential of the left-regular
representation $X(\hat\phi) = -\ell _*(X)\hat \phi$.
Hence, if $\{e_\mu\}$ is an orthonormal basis in $\pp$, the rough
Laplacian $\bar\nabla^*\bar\nabla$
translates into the sum $- \ell _*(e_\mu)\ell _*(e_\mu)   = (-\Cas^G_\ell 
+ \Cas^K_\ell )$. 
Since $\bar\Delta = \bar\nabla^*\bar\nabla + q(\bar R)$ it remains to 
show that $q(\bar R) = -\Cas^K_\ell =  -\Cas^K_\pi$ in order to complete
the proof of the lemma. 

We claim that the differential $i_*:\kk\rightarrow \so(\pp)\cong \Lambda^2\pp$ 
of the isotropy representation $i : K \rightarrow \SO(\pp)$ is given
by $i_*(A) = -\tfrac12 e_\mu \wedge [A,e_\mu]$ for any $A\in
\kk$. Indeed 
$$ 
(\tfrac12 e_\mu \wedge [A,e_\mu])_* X = 
-\tfrac12 B(e_\mu, X)[A,e_\mu] +\tfrac12 B([A,e_\mu], X)e_\mu
= -[A,X]  .
$$
Next we recall that for $X,Y \in \pp $ the curvature $\bar R_{X,Y}$ of the 
canonical connection acts by $-\pi_*([X,Y]^\kk)$ on every associated vector bundle
$EM$, defined by the representation~$\pi$. Hence the curvature operator $\bar{ {R}}$
can be written for any $X,Y\in \pp$ as 
$$
\bar{ {R}} (X\wedge Y) = \tfrac12 e_\mu \wedge \bar R_{X,Y}e_\mu
=-\tfrac12e_\mu \wedge [[X,Y]^\kk,e_\mu]
= i_*([X,Y]^\kk)  .
$$
Let $P_{\SO(\pp)}= G \times_i \SO(\pp)$ be the bundle of orthonormal frames of
$M=G/K$. Then any $\SO(\pp)$-representation $\tilde \pi$ defines a $K$-representation 
by $\pi = \tilde \pi \circ i$. Moreover any 
vector bundle $EM$  associated to $P_{\SO(\pp)}$ via $\tilde \pi$ can
be written as a vector bundle associated via $\pi$ to the
$K$-principle  
bundle $G \rightarrow G/K$, i.e.
$$
EM = P_{\SO(\pp)} \times_{\tilde \pi}E = G \times_\pi E
$$
Let $\{f_\a\}$ be an orthonormal
basis of $\kk$. Then by the definition of $q(\bar R)$ we have
$$
\begin{array}{rcl}
q(\bar R) &=&  \tfrac12	\tilde\pi_* (e_\mu \wedge e_\nu)\,\tilde\pi_* (\bar{R}(e_\mu \wedge e_\nu))
	  =	 \tfrac12  \tilde\pi_* (e_\mu \wedge e_\nu)\,\pi_* ([e_\mu, e_\nu]^\kk)\\[1.5ex]
	 &=&  -\tfrac12	 B([e_\mu, e_\nu], f_\a) \tilde\pi_* (e_\mu \wedge e_\nu)\,\pi_* (f_\a) 
	  =	 -\tfrac12	B(e_\nu, [f_\a, e_\mu]) \tilde\pi_* (e_\mu \wedge e_\nu)\,\pi_* (f_\a) \\[1.5ex]
	 &=&  \tfrac12	\tilde\pi_* (e_\mu \wedge [f_\a, e_\mu])\,\pi_* (f_\a) 
	  =	 -\pi_* (f_\a)\,\pi_* (f_\a)\\[1.5ex]
	 &=&   -\Cas^K_\pi.	
\end{array}
$$
We have shown that $q(\bar R)\in \End(EM)$ acts fibre-wise as
$-\Cas^K_\pi$. Let $Z\in \kk$ and $f\in C^\infty(G;E)^K$, then the
$K$-invariance of $f$ implies $\pi_*(Z)f = -Z(f) = \ell _*(Z)f$ and also
$\Cas^K_\pi = \Cas^K_\ell $, which concludes the proof of the lemma.
\r

It follows from this lemma that the spectrum of $\bar\Delta$ on
sections of $EM$ 
is the set of numbers $\lambda_\g =\|\rho+\g\|^2-\|\rho\|^2$, where $\g$ is the
highest weight of an irreducible $G$-representation $V_\g$ such that 
$\Hom_K(V_\g, E) \neq 0$, i.e. such that the decomposition of $V_\g$,
considered  
as $K$-representation, contains components of the $K$-representation $E$.

\subsection{Nearly K\"ahler deformations and Laplace eigenvalues}
Let $(M,g,J)$ be a  compact simply
connected 6-dimensional nearly K\"ahler manifold not isometric to the
round sphere, with scalar curvature normalized to $\scal_g=30$.
Recall the following result from \cite{dnk}:
\begin{ath}\label{infd} The Laplace operator
$\Delta$ coincides with the Hermitian Laplace operator $\bar\Delta$ on
co-closed primitive $(1,1)$-forms. The space $\NK$ of infinitesimal
deformations of the nearly K\"ahler structure of $M$ is
isomorphic to the eigenspace for the 
eigenvalue $12$ of the restriction of $\Delta$ (or $\bar\Delta$) to
the space of {\em co-closed} primitive $(1,1)$-forms on $M$.
\end{ath}

Assume from now on that $M$ is a 6-dimensional naturally reductive
$3$-symmetric space $G/K$ in the list of 
Gray and Wolf, i.e. $\SU_2\times\SU_2\times \SU_2/\SU_2$, $\SO_5/\U_2$
or $\,SU_3/T^2\,$.  
As was noticed before, the  canonical homogeneous and the canonical Hermitian
connection coincide, since for the later can be shown that is torsion
and its curvature are parallel, a property, which by the
Ambrose-Singer-Theorem characterizes the canonical homogeneous
connection (c.f. \cite{jbb}).  In order to determine the space $\NK$
on $M$ we thus need to apply the previous calculations to
compute the $\bar\Delta$-eigenspace for the eigenvalue $12$ on
primitive $(1,1)$-forms and decide which of these
eigenforms are co-closed.

According to Lemma~\ref{qr} and the decomposition (\ref{pw}) we have
to carry out three steps: 
first to determine the $K$-representation $\Lambda^{1,1}_0 \pp$
defining the bundle  
$\Lambda^{1,1}_0 TM$, then to compute the Casimir eigenvalues with the
Freudenthal formula, 
which gives all possible $\bar\Delta$-eigenvalues and finally to check
whether the 
$G$-representation $V_\gamma$ realizing the eigenvalue $12$ satisfies
$\Hom_K(V_\gamma, \Lambda^{1,1}_0\pp) \neq\{0\}$ and thus really appears
as eigenspace.

Before going on, we make the following useful observation

\begin{elem}\label{spectrum}
Let $(G/K, g)$ be a $6$-dimensional homogeneous strict nearly K\"ahler
manifold of   
scalar curvature $\scal_g = 30$. Then the homogeneous metric $g$ is induced 
from $-\tfrac{1}{12} B$, where $B$ is the Killing form of $G$.
\end{elem}
\begin{proof} 
Let $G/K$ be a 6-dimensional homogeneous strict nearly K\"ahler manifold. Then 
the metric is induced from a multiple of the Killing form, i.e. $G/K$ is a 
normal homogeneous space with $\Ad(K)$-invariant decomposition $\gg =
\kk \oplus \pp$.  
The scalar curvature of the metric $h$ induced by $-B$
may be computed as (c.f. \cite{besse})
$$
\scal_h = \tfrac32 - 3\Cas^K_\lambda
$$
where $\lambda: K \rightarrow \so(\pp)$ is the isotropy
representation. 
From Lemma~\ref{qr} we know that $\Cas^K_\lambda =
-q(\bar R)$, which on the tangent bundle was computed in Lemma~\ref{Rbar} as
$ q(\bar R) = \tfrac{2\scal_h}{15}\, \Id $. Hence we obtain the equation
$\scal_h = \tfrac32 + \tfrac25 \scal_h$ and it follows $\scal_h =
\tfrac52$, i.e. 
the metric $g$ corresponding to $-\tfrac{1}{12}B$ has scalar curvature 
$\scal_g = 30$.
\r


\subsection{The $\bar\Delta$-spectrum on $S^3\times S^3$}

Let $K=\SU_2$ with Lie algebra $\kk = \su_2$ and $G=K\times K \times K
$ with Lie algebra 
$\gg =	\kk\oplus \kk \oplus \kk$. We consider 
the $6$-dimensional manifold $M = G/K$, where $K$ is diagonally
embedded. The tangent space 
at $o=eK$ can be identified with
$$
\pp = \{(X,Y,Z)\in \kk \oplus \kk \oplus \kk \,|\, X+Y+Z=0\}  .
$$
Let $B$ be the Killing form of $\kk$ and define $B_0 = - \tfrac{1}{12}B $. Then
it follows from Lemma~\ref{spectrum} that the invariant scalar product
$$ 
B_0((X,Y,Z),(X,Y,Z)) = B_0(X,X) + B_0(Y,Y) + B_0(Z,Z)
$$
defines a normal metric, which is the homogeneous nearly K\"ahler
metric $g$ of scalar  
curvature $\scal_g = 30$.

The canonical almost complex structure on the $3$-symmetric space $M$,
corresponding to the 3rd order $G$-automorphism $\sigma$, with 
$\sigma(k_1,k_2,k_3) =(k_2,k_3,k_1)$, is defined as
$$
J(X,Y,Z) = \tfrac{2}{\sqrt 3} (Z,X,Y) + \tfrac{1}{\sqrt{3}}(X,Y,Z)  .
$$
The $(1,0)$-subspace $\pp^{1,0}$ of $\pp^\CM$
defined by $J$ is 
isomorphic to the complexified adjoint representation of $SU_2$ on
$\su_2^{\CM}$.  
Let $E=\CM^2$ 
denote the standard representation of $\SU_2$ (notice that $E\cong\bar
E$ because every $\SU_2\cong\Sp_1$ representation is quaternionic). 
\begin{elem}
The $\SU_2$-representation defining the bundle $\Lambda^{(1,1)}_0TM$ splits
into the irreducible summands $\Sym^4E$ and $\Sym^2E$.
\end{elem}
\begin{proof}
The defining $\SU_2$-representation of $\Lambda^{(1,1)}TM$ is 
$\pp^{1,0}\otimes \pp^{0,1} \cong \Sym^2E \otimes
\Sym^2E\cong\Sym^4E\oplus\Sym^2E\oplus\Sym^0E$ from the Clebsch-Gordan
formula. Since we are  
interested in primitive $(1,1)$-forms, we still	have to delete the
trivial summand $\Sym^0E\cong\CM$. 
\r
Since $G= \SU_2\times \SU_2 \times \SU_2$, every
irreducible $G$-representation is isomorphic to one of the representations
$V_{a,b,c}= \Sym^aE \otimes \Sym^bE \otimes \Sym^cE$. The Casimir 
operator of the $\SU_2$-representation $\Sym^kE$ (with respect to $B$) 
is $-\tfrac18k(k+2)$ and the Casimir operator of $G$ is the sum of the
three $\SU_2$-Casimir operators. Hence all possible $\bar\Delta$-eigenvalues
with respect to the metric $B_0$ are of the form
\beq\label{spectrum1}
 \tfrac32(a(a+2) + b(b+2) + c(c+2)) .
\eeq
for non-negative integers $a,b,c$. It is easy to check that the eigenvalue
$12$ is obtained  
only for $(a,b,c)$ equal to $(2,0,0),\ (0,2,0)$ or $(0,0,2)$. 
The restrictions to $\SU_2$ (diagonally embedded in
$G$) of the three corresponding $G$-representations are all equal 
to the $\SU_2$-representation $\Sym^2 E$, thus $\dim \Hom_{\SU_2}(V_{2,0,0},
\Lambda^{(1,1)}_0\pp ) = 1$, and similarly for the two other summands.
Hence the eigenspace of $\bar\Delta$ on primitive $(1,1)$-forms for the 
eigenvalue 12 is isomorphic	 to	 $V_{2,0,0}\oplus
V_{0,2,0}\oplus V_{0,0,2}$ 
and its dimension, i.e. the multiplicity of the eigenvalue $12$, is equal to
$9$.

Since the isometry group of the nearly K\"ahler manifold
$M = \SU_2\times \SU_2 \times \SU_2/ \SU_2$ has dimension $9$, the
inequality (\ref{ineq}) yields 
$$\dim(\NK)\le
\dim(\Omega^{(1,1)}_0(12))-\dim(i(M))-\dim(\Omega^0(12))=
-\dim(\Omega^0(12))\le0.$$   
We thus have obtained the following

\begin{ath} The homogeneous nearly K\"ahler structure on $S^3\times
  S^3$ does not admit any infinitesimal nearly K\"ahler deformations.
\end{ath}

Finally we remark that there are also no infinitesimal Einstein
deformations neither.
In \cite{enk} we showed that the space of infinitesimal Einstein deformations
of a nearly K\"ahler metric $g$, with normalized scalar curvature  $\scal_g
=30$, is isomorphic to the direct sum of 
$\bar\Delta$-eigenspaces 
of primitive co-closed $(1,1)$-forms for the eigenvalues $2, 6$ and
$12$.  It is clear from 
(\ref{spectrum1}) that neither $2$ nor $6$ can be realized as
$\bar\Delta$-eigenvalues. 

\begin{ecor}
 The homogeneous nearly K\"ahler metric on $S^3\times
  S^3$ does not admit any infinitesimal Einstein deformations.
\end{ecor}


\subsection{The $\bar\Delta$-spectrum on $\CM P^3$}

In this section we consider the complex projective space $\CM P^3 =
\SO_5/\U_2$, 
where $\U_2$ is embedded by $\U_2\subset \SO_4\subset \SO_5$. Let
$G=\SO_5$ with 
Lie algebra $\gg$ and $K=\U_2$ with Lie algebra $\kk$. We denote the
Killing form 
of $G$ with $B$. Then we have the $B$-orthogonal decomposition $\gg
=\kk \oplus \pp$, 
where $\pp$ can be identified with the tangent space in $o=eK$. The space $\pp$
splits as $\pp=\mm\oplus\nn$, where $\mm$ resp. $\nn$ can be
identified with the  
horizontal resp. vertical tangent space at $o$ of the twistor space
fibration $\SO_5/\U_2\rightarrow \SO_5/\SO_4 = S^4$. We know from
Lemma~\ref{spectrum} that $B_0 =-\tfrac{1}{12}B$ defines the homogeneous nearly 
K\"ahler metric $g$ of scalar curvature $\scal_g = 30$.

Let $\{\e_1, \e_2\}$ denote the canonical basis of $\RM^2$. Then the positive roots 
of $\SO_5$ are
$\a_1=\e_1,\,\a_2=\e_2,\,\a_3= \e_1+\e_2, \,\a_4=\e_1-\e_2$,\; with
$\;\rho = \tfrac32\e_1+\tfrac12\e_2$. 
Let $\gg^\a\subset \gg^\CM$ be the root space corresponding to the
root $\alpha$. Then 
$$
\mm^\CM = \gg^{\a_1} \oplus \gg^{-\a_1} \oplus \gg^{\a_2} \oplus
\gg^{-\a_2}, \qquad 
\nn^\CM =  \gg^{\a_3} \oplus \gg^{-\a_3}.
$$

The invariant almost complex structure $J$ may be defined by
specifying the $(1,0)$-subspace  
$\pp^{1,0}$ of $\pp^\CM$:
$$
\pp^{1,0}=\{X-iJX\ |\ X\in\pp\}=\gg^{\a_{1}}\oplus \gg^{\a_{2}}\oplus
\gg^{-\a_{3}}  ,
$$
It follows that $J$ is not integrable, since the restricted root
system $\{ \a_{1},\a_{2},-\a_{3}\}$ is not closed under addition
(cf. \cite{HB}). We note that replacing $-\a_3$ by $\a_3$ yields an
integrable almost complex structure. This corresponds to the well-known fact
that on the twistor space the non integrable almost complex structure $J$  is
transformed into the the integrable one by replacing $J$ with $-J$ on the
vertical tangent space.

Let $\CM_k$ denote the $\U_1$-representation on $\CM$ defined by 
$(z,v)\mapsto z^k v$, for $v\in \CM$ and $z \in \U_1 \cong \CM^*$. Then,
since $\U_2 = (\SU_2\times \U_1)/\ZM_2$, any irreducible $\U_2$-representation
is of the form $E_{a,b}= \Sym^aE \otimes \CM_b$, with $a \in \NM,\ b \in \ZM$
and $a \equiv b \mod 2$. As usual let $E = \CM^2$ denote the standard
representation 
of $\SU_2$.

With this notation we obtain the following decomposition of $\pp^{1,0}$ considered as
$\U_2$-representation
\beq\label{p10}
\pp^{1,0} \cong E_{0,-2} \oplus E_{1,1} 
\qquad\mbox{with}\qquad
E_{0,-2} \cong \gg^{-\a_{3}}
\quad\mbox{and}\quad
 E_{1,1}  \cong \gg^{\a_{1}}\oplus \gg^{\a_{2}}
.
\eeq
Since $\pp^{0,1}$ is obtained from $\pp^{1,0}$ by conjugation we have 
$\pp^{0,1} \cong E_{0,2} \oplus E_{1,-1}$. The defining $\U_2$-representation
of $\Lambda^{(1,1)}TM$ is $\pp^{1,0}\otimes \pp^{0,1}$, which obviously 
decomposes into $5$ irreducible summands, among which, two are
isomorphic to the trivial representation
$E_{0,0}$. Considering only primitive $(1,1)$-forms we still have to delete
one of the trivial summands and obtain

\begin{elem}
The $\U_2$-representation defining the bundle $\Lambda^{(1,1)}_0TM$ has
the following decomposition into irreducible summands
$$
\Lambda^{(1,1)}_0\pp =	E_{0,0} \oplus	 E_{1,3} \oplus	 E_{1,-3}
\oplus  E_{2,0}.  
$$
\end{elem}

Let $V_{a,b}$ be an irreducible $\SO_5$-representation of highest weight 
$\g = (a,b)$ with $a,b \in \NM$ and $a \ge b \ge 0$, e.g.
$V_{1,0}=\Lambda^1$ and $V_{1,1}=\Lambda^2$.
The scalar product induced by the 
Killing form $B$ on the dual
$\tt^* \cong \RM^2$ of the maximal torus of $\SO_5$ is $-\tfrac16$
times the Euclidean  
scalar product. By the Freudenthal
formula we thus get
\beq\label{c2}
\Cas_{V_{a,b}}=\la\g,\g+2\rho\ra_B
=-\tfrac{1}{6}(a(a+3)+b(b+1)).
\eeq
Notice that we have $V_{1,1}=\so_5^\CM$ and $\Cas_{V_{1,1}}=-1$, which
is consistent with $\Cas^G_\ad = -1$.

It follows (c.f. Remark \ref{rem1}) that all possible
$\bar\Delta$-eigenvalues with respect to 
the metric induced by $B_0$
are of the form $2(a(a+3)+b(b+1))$. The eigenvalue
$12$ is realized 
if and only if $(a,b)=(1,1)$. We still have to decide whether the
$\SO_5$-representation $V_{1,1}$ 
actually appears in the decomposition (\ref{pw}) of
$L^2(\Lambda^{1,1}_0TM)$. However this 
follows from

\begin{elem}\label{restriction}
The $\SO_5$-representation $V_{1,1}$ restricted to $\U_2\subset \SO_5$ has the
following decomposition as $\U_2$-representation:
$$
V_{1,1} 
\cong 
(E_{0,0} \oplus	   E_{2,0}) \oplus
(E_{0,-2} \oplus E_{1,1} \oplus E_{0,2} \oplus	E_{1,-1})
$$
and in particular
$$
\dim \Hom_{\U_2}(V_{1,1},\Lambda^{1,1}_0\pp^\CM) = 2 
\qquad\mbox{and}\qquad
\dim \Hom_{\U_2}(V_{1,1},\CM) =	 1
.
$$
\end{elem}
\begin{proof}
We know already that $V_{1,1} = \so_5^\CM$ is the complexified adjoint 
representation and that $\so_5^\CM = \uu_2^\CM \oplus (\pp^{1,0}\oplus
\pp^{0,1})$. 
The decomposition of the last two summands is contained in
(\ref{p10}). Hence it 
remains to explicit the adjoint representation of $\U_2$ on $\uu_2^\CM$. It is
clear that its restriction to $\U_1$ acts trivially, whereas
its restriction to $\SU_2$ decomposes into $\CM \oplus
\su_2^\CM$, i.e. $\uu^\CM_2 
\cong E_{0,0} \oplus	  E_{2,0}$. 
\r

The eigenspace of $\bar\Delta$ on primitive $(1,1)$-forms for the
eigenvalue $12$ is thus isomorphic to the sum of two copies of $V_{1,1}$, i.e.
the eigenvalue $12$ has multiplicity $2\cdot 10 = 20$.

It is now easy to calculate the smallest eigenvalue and the corresponding 
eigenspace of the Laplace operator $\Delta$ on non-constant functions.
We do this for $\bar\Delta$, which coincides with $\Delta$ on functions.
Then we have to replace	 $\Lambda^{(1,1)}_0\pp$ in the calculations
above with the trivial representation $\CM$ and to look for 
$\SO_5$-representations $V_{a,b}$ containing the zero weight. 
It follows from Lemma~\ref{restriction} and (\ref{c2}) that the 
$\Delta$-eigenspace on functions $\Omega^0(12)$ is isomorphic to 
$V_{1,1}$ and is thus $10$-dimensional. 
Since the dimension of the isometry group of the nearly K\"ahler manifold
$\SO(5)/\U_2$ is $10$, the inequality (\ref{ineq}) shows
that 
$$\dim(\NK)\le\dim(\L^{(1,1)}_0(12))-\dim(i(M))-\dim(\Lambda^0(12))=
20-10-10=0,$$  
so there are no infinitesimal nearly K\"ahler deformations in this
case neither.

Finally, we remark like before that there are also no other
infinitesimal Einstein 
deformations, since by (\ref{c2}), the eigenvalues $2$ and
$6$ do not occur in the spectrum of $\bar\Delta$ on
$\L^{(1,1)}_0M$. Summarizing, we have obtained the following:

\begin{ath} The homogeneous nearly K\"ahler structure on
  $\CP^3=\SO_5/\U_2$ does not admit any infinitesimal nearly K\"ahler
  or Einstein deformations. 
\end{ath}


\subsection{The $\bar\Delta$-spectrum on the flag manifold $F(1,2)$}

In this section we consider the flag manifold $M= SU_3/T^2$, where
$T^2\subset \SU_3$ is the maximal torus.
Let $\gg=\su_3$ and let $\kk=\tt$, the Lie algebra of $ T^2$. 
We have the decomposition
$$
\gg=\kk\oplus\pp\qquad\hbox{and}\qquad \pp=\mm\oplus\nn.
$$
Denoting by  $E_{ij}$, $S_{ij}$ are "real and imaginary" part of the projection
of the vector $X_{ij}\in\gl_3$ (equal to 1 on $i$th row and $j$th
column and 0 elsewhere) onto $\su_3$:
$$
E_{ij}=X_{ij}-X_{ji}\qquad S_{ij}=i(X_{ij}+X_{ji}),
$$
the subspaces $\mm$ and $\nn$ are explicitly given by
$$
\mm=\span\{E_{12}, S_{12}, E_{13},
S_{13}\}=\span\{e_1,e_2,e_3,e_4\},
$$
$$
\nn=\span\{E_{23},S_{23}\}=\span\{e_5,e_6\}.
$$
The dual of the Lie algebra $\tt$ of the maximal torus $T^2$ can
be identified with 
$$
\tt^*\cong\{(\l_1,\l_2,\l_{3})\in\RM^{3}\ |\ \l_1+\l_2+\l_3=0\}.
$$
If $\{\e_i\}$ denotes the canonical basis in $\RM^{3}$ then the set of
positive roots is given as
$\; \phi ^+=\{\a_{ij}=\e_i-\e_j\ |\ 1\le i<j\le 3\}$
and the half-sum of the positive roots is $\rho = \e_1 -\e_3$

Let $B$ denote the Killing form of $\SU_3$. By Lemma \ref{spectrum},
$B_0 = -\tfrac{1}{12}B $ 
defines the homogeneous nearly K\"ahler metric $g$ of scalar curvature
$\scal_g = 30$.

The almost complex structure $J$ is explicitly defined on $\pp$ by
$$
J(e_1)=e_2,\qquad J(e_3)=-e_4,\qquad J(e_5)=e_6.
$$
Alternatively we may define the $(1,0)$-subspace of $\pp^\CM$:
$$
\pp^{1,0}=\gg^{\a_{12}}\oplus \gg^{\a_{31}}\oplus
\gg^{\a_{23}}=\span\{X_{12},X_{31},X_{23}\},
$$
where $\gg^\a$ is the root space for $\a$.
It follows that $J$ is not integrable, since the restricted root
system $\{ \a_{12},\a_{31},\a_{23}\}$ is not closed under addition
(c.f. \cite{HB}).

Let $E=\CM^{3}$ be the standard representation of $\SU_{3}$ with
conjugate representation $\bar E$.
Any irreducible representations of $\SU_3$ is isomorphic to one of
the representations
$$
V_{k,l}:=(\Sym^k E\otimes\Sym^l\bar E)_0,
$$
where the right hand side denotes the kernel of the contraction map
$$\Sym^k E\otimes\Sym^l\bar E\to \Sym^{k-1} E\otimes\Sym^{l-1}\bar
E,$$
i.e. $V_{k,l}$ is the Cartan summand in $\Sym^k E\otimes\Sym^l\bar E$.
The weights of $\Sym ^k E$ are 
$$
a\e_1+b\e_2+c\e_3,\qquad \hbox{ with } a,b,c\ge0, \ a+b+c= k.
$$
If $v_1,v_2,v_3$ are the weight vectors of $E$, then these weights
correspond to the weight vectors $v_1^a\.v_2^b\.v_3^c$ in $\Sym
^k E$.
Since the weights of $\Sym ^l \bar E$ are just minus the weights of 
$\Sym ^l E$, we see that the weights of $V_{k,l}$ are 
\beq\label{w}
(a-a')\e_1+(b-b')\e_2+(c-c')\e_3,\quad 
a,b,c,a',b',c'\ge0, \ a+b+c= k,\  a'+b'+c'=l.
\eeq
From the given definition of the almost complex
structure $J$ it is clear that the $T^2$-representation on $\pp^{1,0}$ splits 
in three one-dimensional $T^2$-representations with the weights
$\a_{12},\a_{31},\a_{23}$. 

Since the weights of a tensor product representation are the sums of weights
of each factor and since $\e_1+\e_2+\e_3=0$ on the Lie algebra of the 
maximal torus of $\SU_3$, we immediately obtain 

\begin{ecor} \label{cor} The weights of the $T^2$-representation on
  $\L^{1,1}\pp\cong \pp^{1,0}\otimes\pp^{0,1}$ are 
$$\pm 3\e_1,\, \pm3\e_2,\, \pm3\e_3, \hbox{and $0$}.$$
\end{ecor}

It remains to compute the Casimir operator of the irreducible
$\SU_3$-representations 
$V_{k,l}$. The highest weight of $V_{k,l}$ is $\gamma=k\e_1-l\e_3$ and
$\rho = \e_1 - \e_3$, 
thus 
\beq\label{c1}
\Cas_{V_{k,l}}=\la
\g,\g+2\rho\ra_B = -\tfrac{1}{6}(k(k+2)+l(l+2)). 
\eeq
Here we use again the Freudenthal formula and the fact that the
Killing form $B$ induces $-\tfrac16$ times the Euclidean 
scalar product on $\tt^* \subset \RM^3$ (easy calculation).
Notice that we have $V_{1,1}=\su_3^{\CM}$ and $\Cas_{V_{1,1}}=-1$, which
is consistent with $\Cas^G_\ad = -1$ as in the previous cases.

It follows that all possible $\bar\Delta$-eigenvalues (with respect to
the metric $B_0$) 
are of the form $2(k(k+2)+l(l+2))$. Obviously the eigenvalue $12$ can
only be obtained 
for $k=l=1$. Moreover, the restriction of the $\SU_3$-representation
$V_{1,1}$ contains  
the zero weight space. In fact, from (\ref{w}), the zero weight
appears in $V_{k,l}$ if and only if 
there exist $a,b,c,a',b',c'\ge0, \ a+b+c= k,\  a'+b'+c'=l$ such that
$(a-a')\e_1+(b-b')\e_2+(c-c')\e_3=0$, which is equivalent to $k=l$.
We see that 
$\dim \Hom_{T^2}(V_{1,1},\Lambda^{(1,1)}_0\pp ) = 2\cdot 2 = 4$.

Hence the eigenspace of $\bar\Delta$ on primitive $(1,1)$-forms for the
eigenvalue $12$ is isomorphic to the sum of four copies of $V_{1,1}$, i.e.
the eigenvalue $12$ has multiplicity $4\cdot 8 = 32$.

Computing the the smallest eigenvalue and the corresponding 
eigenspace of the Laplace operator $\Delta$ on non-constant functions
we find $V_{0,0}$ for the eigenvalue
$0$ and $V_{1,1}$ for the eigenvalue $12$. All other possible representations
give a larger eigenvalue. Hence, the $\Delta$-eigenspace on functions 
$\Omega^0(12)$ is isomorphic to two copies of $V_{1,1}$, i.e. the 
eigenvalue $12$ has multiplicity $8\cdot 2 = 16$.

Since the dimension of the isometry group of the nearly K\"ahler manifold 
$\SU_3/T^2$ is $8$, we obtain from (\ref{ineq})
\beq\label{d}\dim(\NK)\le\dim(\Omega^{(1,1)}_0(12))-\dim(i(M))-
\dim(\Omega^0(12))=8.\eeq  
In the next section we will show by an explicit construction that
actually the equality holds, so the flag manifold has an 8-dimensional
space of infinitesimal nearly K\"ahler deformations.

Before describing this construction we note that there are no
infinitesimal  
Einstein deformations other than the nearly K\"ahler deformations. It
follows from (\ref{c1}) that the eigenvalue $2$ does not occur in the
spectrum of $\bar\Delta$ on $\L^{(1,1)}_0M$. The
eigenvalue $6$ could be realized on the $\SU_3$-representations 
$V=V_{1,0}$ or $V=V_{0,1}$. However it is
easy to check that $\Hom_{T^2}(V, \Lambda^{(1,1)}_0\pp) = \{0\}$.

\begin{ecor} Every infinitesimal Einstein deformation of the
  homogeneous nearly K\"ahler metric on
  $F(1,2)=\SU_3/T^2$ is an infinitesimal nearly K\"ahler
  deformation. 
\end{ecor}


\section{The infinitesimal nearly K\"ahler deformations on $\SU_3/T^2$}

In this section we describe by explicit computation the space of
infinitesimal nearly K\"ahler deformations of the flag manifold
$F(1,2)=\SU_3/T^2$.
The Lie algebra $\u_3$ is spanned by $\{h_1,h_2,h_3, e_1,\ldots,e_6\}$,
where
$$h_1=iE_{11},\qquad h_2=iE_{22},\qquad h_3=iE_{33},$$
$$e_1=E_{12}-E_{21},\qquad e_3=E_{13}-E_{31},\qquad e_5=E_{23}-E_{32},$$
$$e_2=i(E_{12}+E_{21}), \qquad e_4=i(E_{13}+E_{31}), \qquad
e_6=i(E_{23}+E_{32}). $$ 
We consider the bi-invariant metric $g$ on $\SU_3$ induced by $-B/12$,
where $B$ denotes the Killing form of $\su_3$. It is easy to check
that $|e_i|^2$=1 and $|h_i-h_j|^2=1$ with respect to $g$. We extend
this metric to $\U_3$ in the obvious way which makes the frame
$\{e_i,\sqrt2 h_j\}$ orthonormal. This defines a metric, also denoted
by $g$, on the manifold $M=F(1,2)$. From now on we
identify vectors and 1-forms using this metric and use the notation
$e_{ij}=e_i\wedge e_j$, etc.

An easy explicit commutator calculation yields the exterior derivative
of the left-invariant 1-forms $e_i$ on $\U_3$:
\beq
\begin{array}{rcl}
de_1 &=& -2e_2\wedge(h_1-h_2)+e_{35}+e_{46},\\
de_2 &=& 2e_1\wedge(h_1-h_2)+e_{45}-e_{36},\\
de_3 &=& 2e_4\wedge(h_3-h_1)-e_{15}+e_{26},\\
de_4 &=& -2e_3\wedge(h_3-h_1)-e_{25}-e_{16},\\
de_5 &=& -2e_6\wedge(h_2-h_3)+e_{13}+e_{24},\\
de_6 &=& 2e_5\wedge(h_2-h_3)+e_{14}-e_{23}.\\
\end{array}\label{ddd}
\eeq
Let $J$ denote the almost complex structure on $M=F(1,2)$ whose K\"ahler form
is $\o=e_{12}-e_{34}+e_{56}$ (It is easy to check that $\o$, which a
priori is a left-invariant 2-form on $\U_3$, projects to $M$ because
$L_{h_i}\o=0$). $J$ induces an orientation on $M$ with volume form
$-e_{123456}$. Let $\psp+i\psm$ denote the associated complex volume
form on $M$ defined by the $ad_{T^3}$-invariant form
$(e_2+iJe_2)\wedge (e_4+iJe_4)\wedge (e_6+iJe_6)$. Explicitly,
$$\psp=e_{136}+e_{246}+e_{235}-e_{145},\qquad
\psm=e_{236}-e_{146}-e_{135}-e_{245}.$$
Using (\ref{ddd}) we readily obtain 
\beq\label{u}d(e_{12})=-d(e_{34})=d(e_{56})=\psp,\eeq 
so
$$d\o=3\psp,\qquad\hbox{and}\qquad d\psm=-2\o ^2.$$
The pair $(g,J)$ thus defines a nearly K\"ahler structure on $M$ (a
fact which we already knew).

We fix now an element $\xi\in\su_3\subset\u_3$, and denote by $X$
the {\em right-invariant} vector field on $\U_3$ defined by $\xi$.
Consider the functions
\beq\label{fc}x_i=g(X,e_i),\qquad v_i=g(X,h_i).\eeq
The functions $v_i$ are projectable to $M$ and clearly
$v_1+v_2+v_3=0$. Let us introduce the vector fields on $\U_3$
$$a_1=x_6e_5-x_5e_6,\qquad a_2=x_3e_4-x_4e_3,\qquad a_3=x_2e_1-x_1e_2.$$
One can check that they project to $M$. Of course, one has
$$Ja_1=x_5e_5+x_6e_6,\qquad Ja_2=x_3e_3+x_4e_4,\qquad Ja_3=x_1e_1+x_2e_2.$$
The commutator relations in $\SU_3$ yield
\beq \label{v} dv_1=a_2-a_3,\qquad
dv_2=a_3-a_1,\qquad dv_3=a_1-a_2.\eeq 
Using (\ref{ddd}) and some straightforward computations we obtain 
\beq
\begin{array}{rcl}
d(Ja_1)&=(-a_1+a_2+a_3)\i\psp+4(v_2-v_3)e_{56},\\
d(Ja_2)&=(a_1-a_2+a_3)\i\psp+4(v_1-v_3)e_{34},\\
d(Ja_3)&=(a_1+a_2-a_3)\i\psp+4(v_1-v_2)e_{12}.\\
\end{array}\label{j}
\eeq
We claim that the 2-form
\beq\label{fc1}\f=v_1e_{56}-v_2e_{34}+v_3e_{12}\eeq
on $M$ is of type (1,1), primitive, co-closed, and satisfies
$\Delta\f=12\f$. The first two assertions are obvious (recall that
$v_1+v_2+v_3=0$). In order to prove that $\f$ is co-closed, it is
enough to prove that $d\f\wedge\o=0$. Using (\ref{u}) and (\ref{v})
we compute: \bea d\f\wedge\o&=&[(a_2-a_3)\wedge e_{56}-
(a_3-a_1)\wedge e_{34}+(a_1-a_2)\wedge e_{12}]\wedge(e_{12}-e_{34}+e_{56})\\
&=&(a_1-a_2)\wedge e_{1256}- (a_3-a_2)\wedge
e_{1234}+(a_1-a_2)\wedge e_{3456}=0.\eea 
Finally, using (\ref{j}), we
get \bea \Delta\f&=& d^* d\f=-*d*[(a_2-a_3)\wedge e_{56}-
(a_3-a_1)\wedge
e_{34}+(a_1-a_2)\wedge e_{12}]\\
&=&-*d[Ja_2\wedge e_{12}+Ja_3\wedge e_{34}+Ja_3\wedge
e_{56}-Ja_1\wedge e_{12}-Ja_1\wedge e_{34}-Ja_2\wedge e_{56}]\\
&=& -*[d(Ja_2)\wedge(e_{12}-e_{56})+d(Ja_3)\wedge(e_{34}+e_{56})-
d(Ja_1)\wedge(e_{12}+e_{34})]\\
&=&
-*[(a_1+a_2+a_3)\i\psp\wedge(e_{12}-e_{56}+e_{34}+e_{56}-e_{12}-e_{34})\\
&&-2(a_2\i\psp)\wedge(e_{12}-e_{56})-2(a_3\i\psp)\wedge(e_{34}+e_{56})+
2(a_1\i\psp)\wedge(e_{12}+e_{34})\\
&&+4(v_1-v_3)e_{34}\wedge(e_{12}-e_{56})+4(v_1-v_2)e_{12}\wedge(e_{34}+e_{56})\\&&-
4(v_2-v_3)e_{56}\wedge(e_{12}+e_{34})]\\
&=&-*[4(2v_1-v_2-v_3)e_{1234}+4(v_1+v_3-2v_2)e_{1256}+4(2v_3-v_1-v_2)e_{3456}]\\
&=&-*[12v_1e_{1234}-12v_2e_{1256}+12v_3e_{3456}]=12\f.
 \eea

Taking into account the inequality (\ref{d}), we deduce at once the following
\begin{ecor}
The space of infinitesimal nearly K\"ahler deformations of the nearly
K\"ahler structure on $F(1,2)$ is isomorphic to the Lie algebra of
$\SU_3$. More precisely, every right-invariant vector field $X$ on
$\SU_3$ defines an element $\f\in\NK$ via the formulas $(\ref{fc})$ and
$(\ref{fc1})$. 
\end{ecor}


\labelsep .5cm

\end{document}